# Selection of variables and dimension reduction in high-dimensional non-parametric regression

## Karine Bertin

*Departamento de Estadística, CIMFAV, Universidad de Valparaíso,*
*Avenida Gran Bretaña 1091, Playa Ancha Valparaíso, Chile,*
*Phone: (56)32508324, Fax: (56)32508322.*
*e-mail:* karine.bertin@uv.cl

## Guillaume Lecué

*CNRS, LATP, CMI (UMR 6632), Université de Provence, Technopôle Château-Gombert,*
*39 rue F. Joliot Curie, 13453 Marseille Cedex 13, France*
*Phone: +33 (0) 4 91 11 35 58, Fax: +33(0) 4 91 11 35 34.*
*e-mail:* lecue@latp.univ-mrs.fr

**Abstract:** We consider a $l_1$-penalization procedure in the non-parametric Gaussian regression model. In many concrete examples, the dimension $d$ of the input variable $X$ is very large (sometimes depending on the number of observations). Estimation of a $\beta$-regular regression function $f$ cannot be faster than the slow rate $n^{-2\beta/(2\beta+d)}$. Hopefully, in some situations, $f$ depends only on a few numbers of the coordinates of $X$. In this paper, we construct two procedures. The first one selects, with high probability, these coordinates. Then, using this subset selection method, we run a local polynomial estimator (on the set of interesting coordinates) to estimate the regression function at the rate $n^{-2\beta/(2\beta+d^*)}$, where $d^*$, the "real" dimension of the problem (exact number of variables whom $f$ depends on), has replaced the dimension $d$ of the design. To achieve this result, we used a $l_1$ penalization method in this non-parametric setup.

## 1. Introduction

We consider the non-parametric Gaussian regression model

$$Y_i = f(X_i) + e_i, \quad i = 1, \ldots, n,$$

where the design variables (or input variables) $X_1, \ldots, X_n$ are $n$ i.i.d. random variables with values in $\mathbb{R}^d$, the noise $e_1, \ldots, e_n$ are $n$ i.i.d. Gaussian random variables with variance $\sigma^2$ independent of the $X_i$'s and $f$ is the unknown regression function. In this paper, we are interested in the pointwise estimation of $f$ at a fixed point $x = (x_1, \ldots, x_d) \in \mathbb{R}^d$. We want to construct some estimation procedures $\widehat{f_n}$ having the smallest pointwise integrated quadratic risk

$$\mathbb{E}(\widehat{f_n}(x) - f(x))^2 \tag{1}$$

using only the set of data $D_n = (Y_i, X_i)_{1 \leq i \leq n}$.

Assuming that the regression function enjoys some regularity properties around $x$ is a classical assumption for this problem. In this work, we assume $f$ to be $\beta$-Hölderian around $x$. We recall that a function $f : \mathbb{R}^d \longmapsto \mathbb{R}$ is $\beta$-*Hölderian at the point $x$* with $\beta > 0$, denoted by $f \in \Sigma(\beta, x)$, when the two following points hold:

- $f$ is $l$-times differentiable in $x$ (where $l = \lfloor \beta \rfloor$ is the largest integer which is strictly smaller than $\beta$),
- there exists $L > 0$ such that for any $t = (t_1, \ldots, t_n) \in \mathcal{B}_\infty(x, 1)$,

$$|f(t) - P_l(f)(t, x)| \leq L\|t - x\|_1^\beta,$$

  where $P_l(f)(\cdot, x)$ is the Taylor polynomial of order $l$ associated with $f$ at the point $x$, $\| \cdot \|_1$ is the $l_1$ norm and $\mathcal{B}_\infty(x, 1)$ is the unit $l_\infty$-ball of center $x$ and radius 1.

When $f$ is only assumed to be in $\Sigma(\beta, x)$, no estimator can converge to $f$ (for the risk given in equation (1)) faster than

$$n^{-2\beta/(2\beta+d)}. \tag{2}$$

This rate can be very slow when the dimension $d$ of the input variable $X$ is large. In many practical problems, the dimension $d$ can depend on the number $n$ of observations in such a way that the rate (2) does not even tend to zero when $n$ tends to infinity. This phenomenon is usually called the curse of dimensionality. Fortunately, in some of these problems the regression function really depends only on a few number of coordinates of the input variables. We formulate this heuristic by the following assumption:

**Assumption 1.** *There exist an integer $d^* \leq d$, a function $g : \mathbb{R}^{d^*} \to \mathbb{R}$ and a subset $J = \{i_1, \ldots, i_{d^*}\} \subset \{1, \ldots, d\}$ of cardinality $d^*$ such that for any $(x_1, \ldots, x_d) \in \mathbb{R}^d$*

$$f(x_1, \ldots, x_d) = g(x_{i_1}, \ldots, x_{i_{d^*}}).$$

Under Assumption 1, the "real" dimension of the problem is not anymore $d$ but $d^*$. Then, we hope that if $f \in \Sigma(\beta, x)$ (which is equivalent to say that $g$ is $\beta$-Hölderian at the point $x$), it would be possible to estimate $f(x)$ at the rate given in equation (2) where $d$ is replaced by $d^*$, leading to a real improvement of the convergence rate when $d^* << d$. Nevertheless, starting from the data $D_n$, it is not clear that detecting the set of interesting coordinates $J$ is an easy task. To select this set, we use a $l_1$ penalization technique. This technique has been mostly used in the parametric setup (cf. Bickel et al. (2008), Zhao and Yu (2006), Meinshausen and Yu (2008) and references therein). In the present work, we adapt it to the non-parametric setup and we obtain our first result in this theorem which is a short version of Theorem 1.

**Theorem A** (selection of the subset *J*). *Under Assumption [1] it is possible to construct, only from the data $D_n$, a subset $\widehat{J} \subset \{1, \ldots, d\}$ such that, with probability greater than $1 - c_0 \exp(c_0 d - c_1 n h^{d+2})$ (for a free parameter $0 < h < 1$),*

$$\widehat{J} = J.$$

Once the set $J$ is empirically determined with high probability, we then run a classical local polynomial estimation procedure on the set of indices $\widehat{J}$ to obtain the following theorem which is a short version of Theorem [2].

**Theorem B** (estimation of *f*). *For any $f \in \Sigma(\beta, x)$, with $\beta > 1$, satisfying Assumption [1], it is possible to construct, only from the data $D_n$, an estimation procedure $\widehat{f}_n$ such that*

$$\mathbb{P}[|\widehat{f}_n(x) - f(x)| \geq \delta] \leq c \exp(-c\delta^2 n^{2\beta/(2\beta+d^*)}), \forall \delta > 0$$

*where c does not depend on n.*

The last theorem proves that it is possible, only from the set of data, to reduce and to detect the "real" dimension of the problem under Assumption [1].

The problem we consider in the paper is called a high-dimensional problem. In the last years, many papers have studied these kinds of problems and summarizing here the state of the art is not possible (we refer the reader to the bibliography of Lafferty and Wasserman (2008)). We just mention some papers. In Bickel and Li (2007); Levina and Bickel (2005); Belkin and Niyogi (2003); Donoho and Grimes (2003), it is assumed that the design variable $X$ belongs to a low dimensional smooth manifold of dimension $d^* < d$. All of these work are based on heuristics techniques. In Lafferty and Wasserman (2008), the same problem as the one considered here is handled. Their strategy is a greedy method that incrementally searches through bandwidth in small steps. If the regression $f$ is in a Sobolev ball of order 2, their procedure is nearly optimal for the pointwise estimation of $f$ in $x$. It achieves the convergence rate $n^{-4/(4+d^*+\epsilon)}$ for every $\epsilon > 0$, when $d = \mathcal{O}(\log n / \log \log n)$ and $d^* = \mathcal{O}(1)$. Our procedure improves this result. First, the optimal rate of convergence is achieved. Second, the regression function does not have to be twice differentiable (actually we have the result for any $\beta > 1$). Third, the dimension $d$ can be taken of the order of $\log n$.

The paper is organized as follows. In the coming section, we construct the procedures announced in Theorem A and B. The exact version of Theorem A and B are gathered in Section [3]. Their proofs are given in Section [4].

## 2. Selection and estimation procedures

Our goal is twofold. First, we want to determine the set of indices $J = \{i_1, \ldots, i_{d^*}\}$. Second, we want to construct an estimator of the value $f(x)$ that converges to the rate $n^{-2\beta/(2\beta+d^*)}$ when $f \in \Sigma(\beta, x)$ for $\beta > 1$. To achieve the first goal, we use a $l_1$ penalization of local polynomial estimators.

### *2.1. Selection procedure*

We consider the following set of vectors

$$\bar{\Theta}(\lambda) = \arg\min_{\theta \in \mathbb{R}^{d+1}} \left[ \frac{1}{nh^d} \sum_{i=1}^{n} \left( Y_i - U\left(\frac{X_i - x}{h}\right)\theta \right)^2 K\left(\frac{X_i - x}{h}\right) + 2\lambda\|\theta\|_1 \right], \tag{3}$$

where $U(v) = (1, v_1, \ldots, v_d)$ for any $v = (v_1, \ldots, v_d)^t \in \mathbb{R}^d$, $\|\theta\|_1 = \sum_{j=0}^{d}|\theta_j|$ for any $\theta = (\theta_0, \ldots, \theta_d)^t \in \mathbb{R}^{d+1}$, $h > 0$ is called the *bandwidth*, $\lambda > 0$ is called the *regularization parameter* and $K : \mathbb{R}^d \longrightarrow \mathbb{R}$ is called the *kernel*. We will explain how to choose the parameters $h$ and $\lambda$ in what follows. In the following, we denote $U_0(v) = 1$ and $U_i(v) = v_i$, for $i = 1, \ldots, d$ for any $v = (v_1, \ldots, v_d) \in \mathbb{R}^d$. The kernel $K$ is taken such that the following set of assumptions holds:

**Assumption 2.** *The kernel $K : \mathbb{R}^d \longrightarrow \mathbb{R}$ is symmetric, supported in $\mathcal{B}_\infty(0, 1)$, the matrix $\left(\int_{\mathbb{R}^d} K(y)U_i(y)U_j(y)dy\right)_{i,j \in \{0,\ldots,d\}}$ is diagonal with positive coefficients independent of $d$ in the diagonal and there exists a constant $M_K \geq 1$ independent of $d$ which upper bounds the quantities $\max_{u \in \mathbb{R}^d} |K(u)|$, $\max_{u \in \mathbb{R}^d} K(u)^2$, $\max_{u \in \mathbb{R}^d} |K(u)|\|u\|_2^2$, $\max_{u \in \mathbb{R}^d} |K(u)|\|u\|_2^2$, $\int_{\mathbb{R}^d} K(y)^2(1+\|y\|_2^2)dy$, $\int_{\mathbb{R}^d} |K(u)|^2 \times \|u\|_1^4 du$ and $\int_{\mathbb{R}^d} K(y)^2(U_i(y)U_j(y))^2 dy$.*

Note that for example the uniform kernel $K(u) = \frac{1}{2^d}\mathbb{1}_{\{u \in B_\infty(0,1)\}}$ satisfies the Assumption 2.

Any statistic $\bar{\theta} \in \bar{\Theta}(\lambda)$ is a $l_1$ penalized version of the classical local polynomial estimator. Usually, for the estimation problem of $f(x)$, only the first coordinate of $\bar{\theta}$ is used. Here, for the selection problem, we will use all the coordinates except the first one. We denote by $\widehat{\theta}$ the vector of $\mathbb{R}^d$ made of the $d$ last coordinates of $\bar{\theta}$.

We expect the vector $\widehat{\theta}$ to be sparse (that is with many zero coordinates) such that the set of all the non-zero coordinates of $\widehat{\theta}$, denoted by $\widehat{J}$, will be the same as the set $J$ of all the non-zero coordinates of $(\theta_1^*, \ldots, \theta_d^*)^t$ where $\theta_i^* = h\partial_i f(x)$, for $i \in \{1, \ldots, d\}$, and $\partial_i f(x)$ stands for the $i$−th derivative of $f$ at point $x$. We remark that, under Assumption 1, the vector $(\theta_1^*, \ldots, \theta_d^*)^t$ is sparse.

Note that, the estimator $\bar{\theta} \in \bar{\Theta}(\lambda)$ may not be unique (depending on $d$ and $n$). Hence, the subset selection method may provide different subsets $\widehat{J}$ depending on the choice of $\bar{\theta}$. Nevertheless, Theorem 2 holds for any subset $\widehat{J}$, whatever is the vector $\bar{\theta}$ chosen in $\bar{\Theta}(\lambda)$.

We also consider another selection procedure close to the previous one which requires less assumption on the regression function. We just need to assume that there exists $f_{\max} > 0$ such that $|f(x)| \leq f_{\max}$. With the same notation, we consider the following set of vectors

$$\bar{\Theta}_2(\lambda) = \arg\min_{\theta \in \mathbb{R}^{d+1}} \left[ \frac{1}{nh^d} \sum_{i=1}^{n} \left( Y_i + f_{max} + Ch - U\left(\frac{X_i - x}{h}\right)\theta \right)^2 K\left(\frac{X_i - x}{h}\right) \right.$$
$$\left. + 2\lambda\|\theta\|_1 \right], \tag{4}$$

where $C$ and $h$ will be given later. We just translate the outputs $Y_i$'s by $f_{max} + Ch$. This translation affects the estimator since the LASSO method is not a linear procedure. We denote by $\widehat{J}_2$, this subset selection procedure.

**Remark 1.** *The $l_1$ penalization technique can be related to the problem of linear aggregation (cf. Nemirovski (2000) and Tsybakov (2003)) in a sparse setup. Indeed, $l_1$ penalization is known to provide sparse estimators if the underlying object to estimate is sparse with respect to a given dictionary. Assumption 1 can be interpreted in terms of sparsity of $f$ w.r.t. to a certain dictionary. For that, we consider the set $\mathcal{F} = \{f_0, f_1, \ldots, f_d\}$ of functions from $\mathbb{R}^d$ to $\mathbb{R}$ where $f_0 = \mathbb{1}$ is the constant function equals to 1 and $f_j(t) = (t_j - x_j)/h$ for any $j \in \{1, \ldots, d\}$ and $t = (t_1, \ldots, t_d) \in \mathbb{R}^d$. The set $\mathcal{F}$ is the dictionary. That is the set within we are looking for the best sparse linear combination of elements in $\mathcal{F}$ approaching $f$ in a neighborhood of $x$. In this setup, the Taylor polynomial of order 1 at point $x$, denoted by $P_1(f)(\cdot, x)$, is a linear combination of the elements in the dictionary $\mathcal{F}$. When $f$ is assumed to belong to $\Sigma(\beta, x)$, the polynomial $P_1(f)(\cdot, x)$ is a good approximation of $f$ in a neighborhood of $x$. Moreover, under Assumption 1, this linear combination is sparse w.r.t. the dictionary $\mathcal{F}$. Thus, we hope that, with high probability, minimizing a localized version of the empirical $L_2$-risk penalized by the $l_1$ norm over the set of all the linear combinations of elements in $\mathcal{F}$ will detect the right locations of the interesting indices $i_1, \ldots, i_{d^*}$ (which correspond to the non-zero coefficients of $P_1(f)(\cdot, x)$ in the dictionary $\mathcal{F}$). That is the main idea behind the procedures introduced in this section since we have:*

$$\bar{\Theta}(\lambda) = \arg\min_{\theta \in \mathbb{R}^{d+1}} \left[ \frac{1}{nh^d} \sum_{i=1}^{n} \left( Y_i - \sum_{j=0}^{d} f_j(X_i)\theta_j \right)^2 K\left( \frac{X_i - x}{h} \right) + 2\lambda\|\theta\|_1 \right],$$

*Of course, we can generalize this approach to other dictionaries (this will lead to other sparsity and regularity properties of $f$) provided that the orthogonality properties of $\mathcal{F}$ (cf. Proposition 1) still hold.*

### 2.2. Estimation procedure

We now construct a classical local polynomial estimator (LPE) (cf. Korostelev and Tsybakov (1993); Tsybakov (1986)) on the set of coordinates $\widehat{J}_2$ previously constructed.

We assume that the selection step is now done. We have at hand a subset $\widehat{J}_2 = \{\widehat{i}_1, \ldots, \widehat{i}_{\widehat{d^*}}\} \subset \{1, \ldots, d\}$ of cardinality $\widehat{d^*}$. For the second step, we consider $\gamma_x$ a polynomial on $\mathbb{R}^{\widehat{d^*}}$ of degree $l = \lfloor \beta \rfloor$ which minimizes

$$\sum_{i=1}^{n} (Y_i - \gamma_x(p(X_i - x)))^2 K^\star\left( p\left(\frac{X_i - x}{h^\star}\right) \right)$$

where $h^\star = n^{-1/(2\beta + \widehat{d^*})}$, $p(v) = (v_{\widehat{i}_1}, \ldots, v_{\widehat{i}_{\widehat{d^*}}})^t$ for any $v = (v_1, \ldots, v_d)^t \in \mathbb{R}^d$ and $K^\star : \mathbb{R}^{\widehat{d^*}} \longrightarrow \mathbb{R}$ is a kernel function. The local polynomial estimator of $f$ at the point $x$ is $\widehat{\gamma}_x(0)$ if $\widehat{\gamma}_x$ is unique and 0 otherwise. We denote by $\widehat{f}(x)$ the projection onto $[-f_{max}; f_{max}]$ of the LPE of $f(x)$. Here, we don't use the other coefficients of $\widehat{\gamma}_x(0)$ like we did in the selection step.

For the estimation step, we use a result on the convergence of multivariate LPE from Audibert and Tsybakov (2007). We recall here the properties of the kernel required in Audibert and Tsybakov (2007) to obtain this result.

**Assumption 3.** *The kernel $K^\star : \mathbb{R}^{\widehat{d^\star}} \longrightarrow \mathbb{R}$ is such that: there exists $c > 0$ satisfying*

$$K^\star(u) \geq c \mathbb{1}_{\|x\|_2 \leq c}, \forall u \in \mathbb{R}^{\widehat{d^\star}}; \int_{\mathbb{R}^{\widehat{d^\star}}} K^\star(u) du = 1;$$

$$\int_{\mathbb{R}^{\widehat{d^\star}}} (1 + \|u\|_2^{4\beta})(K^\star(u))^2 du < \infty; \ \sup_{u \in \mathbb{R}^{\widehat{d^\star}}} (1 + \|u\|_2^{2\beta}) K^\star(u) < \infty.$$

## 3. Results

In this section, we provide the main results of this work. To avoid any technical complexity we will assume that the density function $\mu$ of the design $X$ satisfies the following assumption:

**Assumption 4.** *There exists some constants $\eta$, $\mu_m > 0$, $\mu_M \geq 1$ and $L_\mu > 0$ such that*

- $\mathcal{B}_\infty(x, \eta) \subset \text{supp}(\mu)$ *and $\mu_m \leq \mu(y) \leq \mu_M$ for almost every $y \in \mathcal{B}_\infty(x, \eta)$,*
- *$\mu$ is $L_\mu$-Lipschitzian around $x$, that is for any $t \in \mathcal{B}_\infty(x, 1)$, $|\mu(x) - \mu(t)| \leq L_\mu \|x - t\|_\infty$ (remark that the value $\mu(x)$ is the value of the continuous version of $\mu$ around $x$).*

The first result deals with the statistical properties of the selection procedure. For this step, we require a weaker regularity assumption for the regression function $f$. This assumption is satisfied for any $\beta$-Hlderian function in $x$ with $\beta > 1$.

**Assumption 5.** *There exists an absolute constant $L > 0$ such that the following holds. The regression function $f$ is differentiable and*

$$|f(t) - P_1(f)(t, x)| \leq L\|t - x\|_1^\beta, \quad \forall t \in \mathcal{B}_\infty(x, 1),$$

*where $P_1(f)(\cdot, x)$ is the Taylor polynomial of degree 1 of $f$ at the point $x$.*

To achieve an efficient selection of the interesting coordinates, we have to be able to distinguish the non-zero partial derivatives of $f$ from the null partial derivatives. For that, we consider the following assumption:

**Assumption 6.** *There exists a constant $C \geq 72(\mu_M/\mu_m)LM_K\sqrt{d_0}$ such that $|\partial_j f(x)| \geq C$ for any $j \in J$, where the set $J$ is given in Assumption 1 and $d_0$ is an integer such that $d^\star \leq d_0$.*

**Theorem 1.** *There exists some constants $c_0 > 0$ and $c_1 > 0$ depending only on $L_\mu$, $\mu_m$, $\mu_M$, $M_K$, $L$, $C$ and $\sigma$ for which the following holds. We assume that the regression function $f$ satisfies the regularity Assumption 5, the sparsity Assumption 1 such that the integer $d^*$ is smaller than a known integer $d_0$ and the distinguishable Assumption 6. We assume that the density function $\mu$ of the input variable $X$ satisfies Assumption 4.*

*We consider $\bar{\theta} = (\bar{\theta}_0, \ldots, \bar{\theta}_d) \in \bar{\Theta}(\lambda) \subset \mathbb{R}^{d+1}$ and $\bar{\theta}_2 = ((\bar{\theta}_2)_0, \ldots, (\bar{\theta}_2)_d) \in \bar{\Theta}_2(\lambda) \subset \mathbb{R}^{d+1}$ where $\bar{\Theta}(\lambda)$ and $\bar{\Theta}_2(\lambda)$ are defined in equations (3) and (4) with a kernel satisfying Assumption 2, a bandwidth and a regularization parameter such that*

$$0 < h < \frac{\mu_m}{32(d_0+1)L_\mu M_K} \wedge \eta \ and \ \lambda = 8\sqrt{3M_K}\mu_M L h. \qquad (5)$$

*We denote by $\widehat{J}$ the set $\{j \in \{1, \ldots, d\} : \bar{\theta}_j \neq 0\}$ and by $\widehat{J}_2$ the set $\{j \in \{1, \ldots, d\} : (\bar{\theta}_2)_j \neq 0\}$.*

- *If $|f(x)| > Ch$, where $C$ is defined in Assumption 6 or $f(x) = 0$, then with probability greater than $1 - c_1 \exp(c_1 d - c_0 n h^{d+2})$, $\widehat{J} = J$.*
- *If $|f(x)| \leq f_{max}$, then with probability greater than $1 - c_1 \exp(c_1 d - c_0 n h^{d+2})$, $\widehat{J}_2 = J$.*

We remark that Theorem 1 still holds when we only assume that there exists a subset $J \subset \{1, \ldots, d\}$ such that $\partial_j f(x) = 0$ for any $j \notin J$ instead of the more global Assumption 1.

**Theorem 2.** *We assume that the regression function $f$ belongs to the Hölder class $\Sigma(\beta, x)$ with $\beta > 1$ and satisfies the sparsity Assumption 1 such that the integer $d^*$ is smaller than a known integer $d_0$ and the distinguishable Assumption 6. We assume that the density function $\mu$ of the input variable $X$ satisfies Assumption 4 and $|f(x)| \leq f_{max}$. We assume that the dimension $d$ is such that $d + 2 \leq (\log n)/(-2 \log h)$ ($h$ satisfies (5)).*

*We construct the set $\widehat{J}_2$ of selected coordinates with a kernel, a bandwidth and a regularization parameter as in Theorem 1. The LPE estimator $\widehat{f}(x)$ constructed in subsection 2.2 on the subset $\widehat{J}_2$ and a kernel $K^\star$ satisfying Assumption 3, satisfies*

$$\forall \delta > 0, \mathbb{P}[|\widehat{f}(x) - f(x)| \geq \delta] \leq c_1 \exp\Big(-c_2 n^{\frac{2\beta}{2\beta+d^*}} \delta^2\Big),$$

*where $c_1, c_2 > 0$ are constants independent of $n, d, d^*$.*

Note that, by taking the expectation, we obtain $\mathbb{E}[(\widehat{f}(x) - f(x))^2] \leq c n^{\frac{-2\beta}{2\beta+d^*}}$.

**Remark.** *The selection procedure is efficient provided that $c_1 n h^{d+2} - c_0 d$ tends to infinity when $n$ tends to infinity. Namely, we need (with $0 < h < 1$)*

$$d + 2 < \frac{\log n}{-\log h}. \qquad (6)$$

*It is interesting to note that, for d of the order of $\log n$ (like in (6)), the rate of convergence in (2) does not tend to zero. Therefore, in this case and without any previous selection step, a classical LPE can fail to estimate $f(x)$.*

*A remarkable point of Theorem 1 is that the bandwidth h does not have to tend to 0 when n tends to infinity. This particular behavior does not appear when LPE are used for estimation and not for selection. This can be explained because, we do not need to control any bias term in the selection step. The restriction on h comes only from the fact that we need the dictionary $\mathcal{F}$ to be approximatively orthogonal.*

*Finally, once the set of interesting coordinates is selected, we can use it to run other non-parametric methods to estimate the function f with other pointwise risks or integrated risks and under other smoothness assumptions on f. Note that, by considering other order of the $l_1$-penalized LPE in the selection step, it is easy to find other properties of the function f. For instance, inflection points or convexity of f can be detected with a second order method for the selection step.*

## 4. Proofs

### 4.1. Proof of Theorem 1

First note that, considering only the observations $X_i$ in the neighborhood of $x$, an estimator $\bar{\theta} = (\bar{\theta}_0, \ldots, \bar{\theta}_d) \in \bar{\Theta}(\lambda)$ defined in (3) can be viewed as a Lasso estimator in the linear regression model

$$Z = A\theta^* + \varepsilon, \text{ where } \theta^* = (\theta_0^*, \ldots, \theta_d^*)^t = (f(x), h\partial_1 f(x), \ldots, h\partial_d f(x))^t \quad (7)$$

and, for any $i = 1, \ldots, n$ $\Delta_i := \alpha_i f(X_i) - A_i \theta^*$ and $\alpha_i := \frac{1}{\sqrt{nh^d}} K^{1/2} \left( \frac{X_i - x}{h} \right)$, the output vector $Z$ of $\mathbb{R}^n$ has for coordinates $Z_i := \alpha_i Y_i, i = 1, \ldots, n$, the lines of the design matrix $A \in \mathcal{M}_{n,d+1}$ are $A_i := \alpha_i U \left( \frac{X_i - x}{h} \right), i = 1, \ldots, n$ ($U$ is defined after Equation (3)) and the noise vector $\varepsilon$ has $\varepsilon_i = \alpha_i e_i + \Delta_i$ for coordinates. We remark that the noise is not centered. The "localized" bias term $\Delta := (\Delta_1, \ldots, \Delta_n)^t$ has been added to the noise. With this new notation, we have

$$\bar{\Theta}(\lambda) = \arg\min_{\theta \in \mathbb{R}^{d+1}} \|Z - A\theta\|_2^2 + 2\lambda\|\theta\|_1$$

where $\bar{\Theta}(\lambda)$ has been introduced in equation (3) and $\forall z = (z_1, \ldots, z_n) \in \mathbb{R}^n, \|z\|_2^2 = \sum_{i=1}^n z_i^2$.

For the same reason, an estimator $\bar{\theta}_2 \in \bar{\Theta}_2(\lambda)$ defined in (4) can be viewed as a Lasso estimator in the linear regression model

$$\check{Z} = A\check{\theta}^* + \varepsilon, \text{ where } \check{\theta}^* = (f(x) + f_{max} + Ch, h\partial_1 f(x), \ldots, h\partial_d f(x))^t$$

and $\check{Z}$ has for coordinates $\check{Z}_i = \alpha_i(Y_i + f_{max} + Ch), i = 1, \ldots, n$. Note that the $\Delta_i$'s are not affected by this translation.

We start by studying $\bar{\theta} \in \bar{\Theta}(\lambda)$ when $|f(x)| \geq Ch$ and $\bar{\theta}_2 \in \bar{\Theta}_2(\lambda)$ when $|f(x)| \leq f_{max}$. The study of $\bar{\theta}$ when $f(x) = 0$ will be discussed at the end. Note that, in both the considered cases, we have $|\theta_0^*| \geq Ch$ and $|\tilde{\theta}_0^*| \geq Ch$. This fact will be used in what follows. We first study $\tilde{\theta}$ when $|f(x)| \geq Ch$. The study of $\tilde{\theta}_2$ when $|f(x)| \leq f_{max}$ is the same with the translated data $\tilde{Y}_i = Y_i + f_{max} + Ch$ and $\tilde{f} = f + f_{max} + Ch$. Note that $\tilde{f}$ and $f$ have the same partial derivatives thus $\theta^*$ and $\tilde{\theta}^*$ have the same last $d$ coordinates which are the only ones of interest for the selection step.

Proving Theorem 1 can be viewed as a problem of sign consistency of the Lasso estimator $\hat{\theta} = (\hat{\theta}_1, \ldots, \hat{\theta}_d)$ (the vector made of the $d$ last coordinates of $\bar{\theta}$). To solve this problem, we follow the lines of Zhao and Yu (2006). We remark that, we treat carefully the problem of uniqueness of the LASSO contrary to the work of Zhao and Yu (2006) where uniqueness of the LASSO estimator was assumed.

We first treat the problem of uniqueness of the LASSO. We introduce the function

$$\phi(\theta) := \|Z - A\theta\|_2^2 + 2\lambda \|\theta\|_1, \quad \forall \theta \in \mathbb{R}^{d+1} \tag{8}$$

and we say that $\theta \in \mathbb{R}^{d+1}$ satisfies the system $(S)$ when

$$\forall j = 0, \ldots, d, \begin{cases} (A_{.j})^t (Z - A\theta) = -\lambda \mathrm{sign}(\theta_j) & \text{if} \quad \theta_j \neq 0 \\ |(A_{.j})^t (Z - A\theta)| \leq \lambda & \text{if} \quad \theta_j = 0 \end{cases}$$

where, for any $j \in \{0, \ldots, d\}$, the vector $A_{.j}$ is the $j$-th column of $A$. It is known that $\theta \in \mathbb{R}^{d+1}$ belongs to $\bar{\Theta}(\lambda)$ if and only if $\theta$ satisfies the system $(S)$.

**Lemma 1.** *If $\bar{\theta} \in \mathbb{R}^{d+1}$ and $\bar{\theta}^{(2)} \in \mathbb{R}^{d+1}$ are two solutions of $(S)$ then $A\bar{\theta} = A\bar{\theta}^{(2)}$.*

*Proof of Lemma 1.* We denote by $S(\bar{\theta})$ the set $\{j \in \{0, \ldots, d\} : \bar{\theta}_j \neq 0\}$. For any $v \in \mathbb{R}^{d+1}$, we have

$$\phi(\bar{\theta}+v) - \phi(\bar{\theta}) = 2\lambda \sum_{j \in S(\bar{\theta})} |\bar{\theta}_j + v_j| - |\bar{\theta}_j| - v_j \mathrm{sign}(\bar{\theta}_j) + 2\lambda \sum_{j \notin S(\bar{\theta})} |v_j| - \eta_j v_j + \|Av\|_2^2,$$

where $\eta_j = \lambda^{-1}(A_{.j})^t(Z - A\bar{\theta})$. For any $j \in S(\bar{\theta})$, we have $|\bar{\theta}_j + v_j| - |\bar{\theta}_j| - v_j \mathrm{sign}(\bar{\theta}_j) \geq 0$ and for any $j \notin S(\bar{\theta})$, we have $|\eta_j| \leq 1$ so $|v_j| - \eta_j v_j \geq 0$. Hence,

$$\phi(\bar{\theta} + v) - \phi(\bar{\theta}) \geq \|Av\|_2^2$$

We take $v \in \mathbb{R}^{d+1}$ such that $\bar{\theta}^{(2)} = \bar{\theta} + v$. The vectors $\bar{\theta}^{(2)}$ and $\bar{\theta}$ are both solutions of $(S)$, thus they are minimizers of $\phi$ and so $\phi(\bar{\theta}^{(2)}) = \phi(\bar{\theta})$. Therefore, we have $\|Av\|_2^2 = 0$. $\square$

Next, we prove a result which deals with the identifiability of the model as well as the uniqueness of the LASSO. We introduce the event

$$\Omega_{01} := \left\{ \forall \theta \in \mathbb{R}^{d+1} : \frac{1}{2} \sqrt{\frac{\mu_m}{2}} \|\theta\|_2 \leq \|A\theta\|_2 \leq 2 \sqrt{\frac{3\mu_M}{2}} \|\theta\|_2 \right\}.$$

**Proposition 1.** *There exists two constants $c_0$ and $c_1$ depending only on $\mu_m, \mu_M$ and $M_K$ such that, under Assumption 2 and the first point of Assumption 4 with $0 < h < \eta$, we have*

$$\mathbb{P}[\Omega_{01}] \geq 1 - c_0 \exp(c_0 d - c_1 n h^d).$$

*Proof of Proposition 1.* Let $\theta \in \mathbb{R}^{d+1}$. We have

$$\|A\theta\|_2^2 = \frac{1}{n} \sum_{k=1}^{n} < \bar{Z}_k, \theta >^2 \quad \text{where } \bar{Z}_k := \sqrt{n}\alpha_k U^t \left( \frac{X_k - x}{h} \right).$$

It is easy to see that $| < \bar{Z}_k, \theta >^2 | \leq (2M_K/h^d)\|\theta\|_2^2$ and $\mathbb{V}(< \bar{Z}_k, \theta >^2) \leq (M_K\mu_M/h^d)\|\theta\|_2^4$. Let $0 < \gamma < 1$ be a number that will be chosen wisely latter. Bernstein's inequality yields that, with probability greater than $1 - 2\exp(-9nh^d\gamma^2\mu_M/(16M_K))$,

$$\left| \|A\theta\|_2^2 - \mathbb{E}\|A\theta\|_2^2 \right| \leq (3/2)\gamma\mu_M\|\theta\|_2^2.$$

Moreover, we have $\mathbb{E}\|A\theta\|_2^2 = \int_{\mathbb{R}^d} K(t)(U(t)\theta)^2\mu(x + ht)dt$. To simplify the proof we will suppose that $(\int_{\mathbb{R}^d} K(t)U_i(t)U_j(t)dt)_{0 \leq i,j \leq d} = I_{d+1}$ but the proof still holds when this matrix is diagonal with positive coefficients independent of $d$ as in Assumption 2. Then we obtain that $\mu_m\|\theta\|_2^2 \leq \mathbb{E}\|A\theta\|_2^2 \leq \mu_M\|\theta\|_2^2$. Thus, with probability greater than $1 - 2\exp(-9nh^d\gamma^2\mu_M/(16M_K))$, we have

$$(\mu_m - (3/2)\gamma\mu_M)\|\theta\|_2^2 \leq \|A\theta\|_2^2 \leq (\mu_M + (3/2)\gamma\mu_M)\|\theta\|_2^2. \tag{9}$$

To control the probability measure of $\Omega_{01}$, we need a uniform control over $\theta \in \mathbb{R}^{d+1}$ of $\|A\theta\|_2^2$. For that we use a classical $\epsilon$-net argument. For the sake of completeness, we recall here this argument. Let $\epsilon > 0$ be chosen wisely later and $N_\epsilon$ be an $\epsilon$-net of $\mathcal{S}^d$ (the unit sphere of $\mathbb{R}^{d+1}$) for the $\|\cdot\|_2$-norm. Using an union bound and equation (9), with probability greater than $1 - 2|N_\epsilon|\exp(-9nh^d\gamma^2\mu_M/(16M_K))$, we have

$$(\mu_m - (3/2)\gamma\mu_M) \leq \|A\theta\|_2^2 \leq (\mu_M + (3/2)\gamma\mu_M), \quad \forall \theta \in N_\epsilon. \tag{10}$$

Now, we want to extend the last result to the whole sphere $\mathcal{S}^d$. Let $\theta \in \mathcal{S}^d$. There exists $\theta_0 \in N_\epsilon$ such that $\|\theta - \theta_0\|_2 \leq \epsilon$. If $\theta \neq \theta_0$, there exists $\theta_1 \in N_\epsilon$ which is $\epsilon$-close to $(\theta - \theta_0)/\|\theta - \theta_0\|_2$. Using this argument recursively, we obtain that there exists a sequence $(\delta_j)_{j \geq 0}$ of non-negative numbers such that $\delta_0 = 1$ and $|\delta_j| \leq \epsilon^j, \forall j \geq 1$ and a sequence $(\theta_j)_{j \geq 0}$ of elements in $N_\epsilon$ such that

$$\theta = \sum_{j=0}^{\infty} \delta_j\theta_j.$$

Thus, for any $\theta \in \mathcal{S}^d$, we have

$$\|A\theta\|_2 = \left\| A\left( \sum_{j=0}^{\infty} \delta_j\theta_j \right) \right\|_2 \leq \sum_{j=0}^{\infty} |\delta_j| \|A\theta_j\|_2 \leq \frac{1}{1 - \epsilon} \max_{\theta \in N_\epsilon} \|A\theta\|_2 \tag{11}$$

and

$$\|A\theta\|_2 = \left\| A\left( \sum_{j=0}^{\infty} \delta_j \theta_j \right) \right\|_2 \geq \|A\theta_0\| - \sum_{j=1}^{\infty} |\delta_j| \|A\theta_j\|_2$$

$$\geq \min_{\theta \in N_\epsilon} \|A\theta\|_2 - \frac{\epsilon}{1-\epsilon} \max_{\theta \in N_\epsilon} \|A\theta\|_2. \tag{12}$$

We take $\gamma = \mu_m/(3\mu_M)$ and $\epsilon = (1/4)\sqrt{\mu_m/(3\mu_M)}$. We know that there exists an absolute constant $c > 0$ such that $|N_\epsilon| \leq (c/2)\epsilon^{-d}$. Using this fact and equations (10), (11) and (12), with probability greater than $1 - c \exp(c_0 d - c_1 nh^d)$, we have

$$\frac{1}{2}\sqrt{\frac{\mu_m}{2}} \leq \|A\theta\|_2 \leq 2\sqrt{\frac{3\mu_M}{2}}, \quad \forall \theta \in \mathcal{S}^d,$$

where $c_0 = \log[4\sqrt{(3\mu_M)/\mu_m}] \vee c$ and $c_1 = \mu_m^2/(M_K\mu_M)$. We complete the proof by applying this result to the vector $\theta/\|\theta\|_2$ for any $\theta \in \mathbb{R}^d - \{0\}$ (the result is obvious for $\theta = 0$). □

We introduce the squared matrix of $\mathcal{M}_{d+1}$

$$\Psi^{(n)} := A^t A \text{ and } \Psi_{ij}^{(n)}$$

$$:= \frac{1}{nh^d} \sum_{k=1}^{n} K\left( \frac{X_k - x}{h} \right) U_i \left( \frac{X_k - x}{h} \right) U_j \left( \frac{X_k - x}{h} \right), \quad \forall i,j = 0, \ldots, d,$$

where we recall that $U_0(v) = 1$ and $U_i(v) = v_i, i = 1, \ldots, d$ for any $v \in \mathbb{R}^d$.

To simplify notation and without loss of generality, we will assume, in what follows, that the interesting indexes are given by the first $d^*$ coordinates. Namely, we will assume (but we did not use it to construct our procedures) that $(i_1, \ldots, i_{d^*}) = (1, \ldots, d^*)$ and then $J = \{1, \ldots, d^*\}$.

We introduce some notation. The vector $\theta^*$ and the matrices $\Psi^{(n)} = A^t A$ and $A$ can be written as

$$\Psi^{(n)} := \begin{pmatrix} \Psi_{11} & \Psi_{12} \\ \Psi_{21} & \Psi_{22} \end{pmatrix}; \ A := (A_{(1)} A_{(2)}) \text{ and } \theta^* := \begin{pmatrix} \theta_{(1)}^* \\ \theta_{(2)}^* \end{pmatrix},$$

where $\Psi_{11} \in \mathcal{M}_{d^*+1}$, $\Psi_{12} \in \mathcal{M}_{d^*+1,d-d^*}$, $\Psi_{21} \in \mathcal{M}_{d-d^*,d^*+1}$, $\Psi_{22} \in \mathcal{M}_{d-d^*}$, $A_{(1)} \in \mathcal{M}_{n,d^*+1}$, $A_{(2)} \in \mathcal{M}_{n,d-d^*}$, $\theta_{(1)}^* \in \mathbb{R}^{d^*+1}$ and $\theta_{(2)}^* \in \mathbb{R}^{d-d^*}$. We remark that, with the notational simplifications, $\theta_{(2)}^*$ is the null vector of $\mathbb{R}^{d-d^*}$.

**Lemma 2.** *On the event $\Omega_{01}$, the following statements hold:*

- *the LASSO selector exists and is unique,*
- *all the eigenvalues of $\Psi^{(n)}$, $\Psi_{11}$ and $\Psi_{22}$ belong to $[\mu_m/8, 6\mu_M]$,*

*Proof of Lemma 2.* For the first point, we use the convexity of the function $\phi$ (introduced in equation (8)) to obtain the existence of a LASSO selector. By Lemma 1, two LASSO selectors are in the kernel of $A$. On the event $\Omega_{01}$, the kernel of $A$ is $\{0\}$. Thus, there is uniqueness of the LASSO on $\Omega_{01}$.

For the second point, we know that the eigenvalues of $\Psi^{(n)}$ are the square of the singular values of $A$. On the event $\Omega_{01}$, the singular values of $A$ belong to $[(1/2)\sqrt{\mu_m/2}, 2\sqrt{3\mu_M/2}]$. This completes the proof for $\Psi^{(n)}$. Now, let $\lambda > 0$ be an eigenvalue of $\Psi_{11}$ and $v_{(1)} \in \mathbb{R}^{d^*+1}$ be an eigenvector associated with $\lambda$. We denote by $u_{(2)}$ the null vector of $\mathbb{R}^{d-d^*}$. We have

$$\lambda \|u_{(1)}\|_2^2 = \|A_{(1)}u_{(1)}\|_2^2 = \|A(u_{(1)}^t u_{(2)}^t)^t\|_2^2 \leq 6\mu_M \|(u_{(1)}^t u_{(2)}^t)^t\|_2^2 = 6\mu_M \|u_{(1)}\|_2^2,$$

thus $\lambda \leq 6\mu_M$. For the same reason, we have $\lambda \geq \mu_m/8$. The proof for $\Psi_{22}$ follows the same argument. □

We consider the event

$$\Omega_{02} := \left\{ \forall j \in \{d^*+1, \ldots, d\}, \forall k \in \{0, \ldots, d^*\} : |(\Psi_{21})_{jk}| \leq 2hL_\mu M_K \right\}. \quad (13)$$

**Lemma 3.** *We assume that Assumption 2 and Assumption 4 hold. There exists a constant $c_3$ depending only on $L_\mu, M_K$ and $\mu_M$ such that the following holds. We have*

$$\mathbb{P}[\Omega_{02}] \geq 1 - 2(d^*+1)(d-d^*)\exp(-c_3 n h^{d+2}).$$

*We take $h$ such that $0 < h < \mu_m[32(d^*+1)L_\mu M_K]^{-1} \wedge \eta$. On the event $\Omega_{01} \cap \Omega_{02}$, we have*

$$\forall j = d^*+1, \ldots, d \quad |(\Psi_{21}(\Psi_{11})^{-1}\overrightarrow{\mathrm{sign}}(\theta^*_{(1)}))_j| < 1/2.$$

*Proof of Lemma 3.* The first point is a direct application of Bernstein's inequality and of the union bound. We use both assumptions of the lemma to upper bound the expectation $|\mathbb{E}(\Psi_{12})_{jk}| \leq hL_\mu M_K$.

For the second part of the lemma, let $j \in \{d^*+1, \ldots, d\}$. On the event $\Omega_{01}$, the maximal eigenvalue of $\Psi_{11}^{-1}$ is smaller than $8/\mu_m$, thus, we have

$$|(\Psi_{21}\Psi_{11}^{-1}\overrightarrow{\mathrm{sign}}(\theta^*_{(1)}))_j| = \Big| \sum_{k=0}^{d^*} (\Psi_{21})_{jk}(\Psi_{11}^{-1}\overrightarrow{\mathrm{sign}}(\theta^*_{(1)}))_k \Big|$$

$$\leq \Big( \sum_{k=0}^{d^*} (\Psi_{21})_{jk}^2 \Big)^{1/2} \|\Psi_{11}^{-1}\overrightarrow{\mathrm{sign}}(\theta^*_{(1)})\|_2$$

$$\leq \sqrt{d^*+1}(2hL_\mu M_K)\frac{8\sqrt{d^*+1}}{\mu_m} < 1/2.$$

□

**Remark.** *If we have $\mathbb{E}\psi^{(n)} = I_{d+1}$, then we don't need any restriction on $h$. Because, in this case, we can obtain that, with high probability, $\forall \theta \in \mathbb{R}^{d+1}$, $(1-\gamma)\|\theta\|_2 \leq \|A\theta\|_2 \leq (1+\gamma)\|\theta\|_2$. Thus, with the same probability, we have*

$$\forall \theta \in \mathbb{R}^{d+1}, (1-\gamma)\|\theta\|_2 \leq \|\Psi\theta\|_2 \leq (1+\gamma)\|\theta\|_2.$$

*By applying the last inequality to the vector $\theta_0 = (\Psi_{11}^{-1}\overrightarrow{\text{sign}}(\theta_{(1)}^*)^t 0^t)^t$, we obtain:*

$$d^* + 1 + \|\Psi_{21}\Psi_{11}^{-1}\overrightarrow{\text{sign}}(\theta_{(1)}^*)\|_2^2 = \|(\overrightarrow{\text{sign}}(\theta_{(1)}^*)^t, [\Psi_{21}\Psi_{11}^t\overrightarrow{\text{sign}}(\theta_{(1)}^*)]^{-1})\|_2^2$$

$$= \|\Psi\theta_0\|_2^2 \leq (1+\gamma)^2\|\theta_0\|_2^2 = (1+\gamma)^2\|\Psi_{11}^{-1}\overrightarrow{\text{sign}}(\theta_{(1)}^*)\|_2^2 \leq (1+\gamma)(d^*+1).$$

*Thus, we get $\|\Psi_{21}\Psi_{11}^{-1}\overrightarrow{\text{sign}}(\theta_{(1)}^*)\|_2^2 \leq \gamma(d^*+1)$. Thus, for $\gamma$ small enough we have $|(\Psi_{21}\Psi_{11}^{-1}\overrightarrow{\text{sign}}(\theta_{(1)}^*))_j| \leq 1/2$. The problem is that, in general, the dictionary cannot satisfies $\mathbb{E}\psi^{(n)} = I_{d+1}$. We just have*

$$(\mu(x) - m_0 h)I_{d+1} \leq \mathbb{E}\Psi^{(n)} \leq (\mu(x) + m_1 h)I_{d+1},$$

*with $m_0$ and $m_1$ two positive constants.*

We consider the following events:

$$\Omega_0 := \Omega_{01} \cap \Omega_{02},$$

$$\Omega_1 := \left\{\forall j = 0, \ldots, d^* : |(\Psi_{11}^{-1}W_{(1)})_j - \lambda(\Psi_{11}^{-1}\overrightarrow{\text{sign}}(\theta_{(1)}^*))_j| < |(\theta_{(1)}^*)_j|\right\}$$

and

$$\Omega_2 := \left\{\forall j = d^*+1, \ldots, d : |(\Psi_{21}\Psi_{11}^{-1}W_{(1)} - W_{(2)})_j| < \lambda/2\right\},$$

where $W_{(1)} = A_{(1)}^t\varepsilon \in \mathbb{R}^{d^*+1}$ and $W_{(2)} = A_{(2)}^t\varepsilon \in \mathbb{R}^{d-d^*}$. For notational simplicity, the indices of the coordinates of any vector in $\mathbb{R}^{d^*+1}$ start from 0 and go to $d^*$, and for any vector in $\mathbb{R}^{d-d^*}$ the indices start from $d^*+1$ and go to $d$. We remark that, we work only on the event $\Omega_0$ on which the minimum eigenvalue of $\Psi_{11}$ is strictly positive. Thus, on this event, $\Psi_{11}$ is regular and so $\Omega_0 \cap \Omega_1$ and $\Omega_0 \cap \Omega_2$ are well defined.

**Proposition 2.** *Let $0 < h < \mu_m[32(d^*+1)L_\mu M_K]^{-1} \wedge \eta$ and Assumption 2 and Assumption 4 hold. The event $\left\{\forall \bar{\theta} \in \mathbb{R}^{d+1} \text{ solution of } (S), \text{ we have } \overrightarrow{\text{sign}}(\widehat{\theta}) = \overrightarrow{\text{sign}}(\theta^*)\right\} \cap \Omega_0$ contains the event $\Omega_0 \cap \Omega_1 \cap \Omega_2$. We recall that for any $\bar{\theta} \in \mathbb{R}^{d+1}$, the vector $\widehat{\theta} = (\bar{\theta}_1, \ldots, \bar{\theta}_d)$ is the vector made of the $d$ last coordinates of $\bar{\theta}$.*

*Proof of Proposition 2.* We consider the linear functional

$$F : \begin{cases} \mathbb{R}^{d^*+1} & \longrightarrow & \mathbb{R}^{d^*+1} \\ \theta & \longmapsto & \theta - \theta_{(1)}^* - \Psi_{11}^{-1}W_{(1)} + \lambda\alpha_{(1)} \end{cases}$$

where we denote by $\alpha_{(1)}$ the vector $\Psi_{11}^{-1}\overrightarrow{\text{sign}}(\theta_{(1)}^*)$. For any $v = (v_0, \ldots, v_d)^t \in \mathbb{R}^{d^*+1}$ and $r = (r_0, \ldots, r_d) \in (\mathbb{R}_+^*)^{d^*+1}$, we set $\mathcal{B}(x,r) = \Pi_{j=0}^{d^*}(x_j - r_j; x_j + r_j)$. For any vector $v = (v_j)_j$, we set $|v| = (|v_j|)_j$. We have $F(\mathcal{B}(\theta_{(1)}^*, |\theta_{(1)}^*|)) = \mathcal{B}(-\Psi_{11}^{-1}W_{(1)} + \lambda\alpha_{(1)}, |\theta_{(1)}^*|)$. On the event $\Omega_1$, we have $0 \in \mathcal{B}(-\Psi_{11}^{-1}W_{(1)} + \lambda\alpha_{(1)}, |\theta_{(1)}^*|)$. Hence, there exists $\bar{\theta}_{(1)} \in \mathcal{B}(\theta_{(1)}^*, |\theta_{(1)}^*|)$ such that $F(\bar{\theta}_{(1)}) = 0$.

That is $\bar{\theta}_{(1)} = \theta^*_{(1)} + \Psi_{11}^{-1}W_{(1)} - \lambda\alpha_{(1)}$ and $|\bar{\theta}_{(1)} - \theta^*_{(1)}| < |\theta^*_{(1)}|$. Thus, we have $\overrightarrow{\text{sign}}(\bar{\theta}_{(1)}) = \overrightarrow{\text{sign}}(\theta^*_{(1)})$ and so

$$\Psi_{11}(\bar{\theta}_{(1)} - \theta^*_{(1)}) = -\lambda\overrightarrow{\text{sign}}(\bar{\theta}_{(1)}). \tag{14}$$

Using Lemma 3, on the event $\Omega_0$, we have for any $j = d^* + 1, \ldots, d$, $|(\Psi_{21}\Psi_{11}^{-1}\overrightarrow{\text{sign}}(\theta^*_{(1)}))_j| < 1/2$. Thus, on the event $\Omega_0 \cap \Omega_2$, we have

$$-\lambda\mathbb{1}_{d-d^*} < \Psi_{21}\Psi_{11}^{-1}W_{(1)} - \lambda\Psi_{21}\alpha_{(1)} - W_{(2)} = \Psi_{21}(\bar{\theta}_{(1)} - \theta^*_{(1)}) - W_{(2)} < \lambda\mathbb{1}_{d-d^*}, \tag{15}$$

where the last inequality is coordinates by coordinates and $\mathbb{1}_{d-d^*}$ is the unit vector of $\mathbb{R}^{d-d^*}$.

We consider the vector $\bar{\theta} = (\bar{\theta}^t_{(1)}, \bar{\theta}^t_{(2)})^t \in \mathbb{R}^{d+1}$ where $\bar{\theta}_{(2)} = \mathbf{0}_{d-d^*}$ is the null vector of $\mathbb{R}^{d-d^*}$. We thus have $\overrightarrow{\text{sign}}(\bar{\theta}) = \overrightarrow{\text{sign}}(\theta^*)$ (because $\overrightarrow{\text{sign}}(\bar{\theta}_{(1)}) = \overrightarrow{\text{sign}}(\theta^*_{(1)})$ and $\bar{\theta}_{(2)} = \mathbf{0}_{d-d^*} = \theta^*_{(2)}$). Moreover equations (14) and (15) are equivalent to say that $\bar{\theta}$ satisfies the system $(S)$.

In particular, we prove that on the event $\Omega_0 \cap \Omega_1 \cap \Omega_2$, there exists $\bar{\theta} \in \mathbb{R}^{d+1}$ solution of $(S)$ such that $\overrightarrow{\text{sign}}(\bar{\theta}) = \overrightarrow{\text{sign}}(\theta^*)$. We complete the proof with the uniqueness of the LASSO on the event $\Omega_0$. $\qquad\square$

**Proposition 3.** *There exists $c_0 > 0$ and $c_1 > 0$ depending only on $L_\mu, \mu_m, \mu_M$, $M_K, L, C$ and $\sigma$ such that the following holds. Under the same assumption as in Theorem 1, we have*

$$\mathbb{P}[\Omega_0 \cap \Omega_1 \cap \Omega_2] \geq 1 - c_1\exp(c_1 d - c_0 nh^{d+2}).$$

*Proof of Proposition 3.* We study the probability measure of $\Omega_2$. We have

$$\Omega_2^c \subset \cup_{j=d^*+1}^d \{|\zeta_j| \geq \lambda/2 - |b_j|\},$$

where $b = (b_{d^*+1}, \ldots, b_d)^t = G(\Delta_1, \ldots, \Delta_n)^t$ and $\zeta = (\zeta_{d^*+1}, \ldots, \zeta_d)^t = G(\alpha_1 e_1, \ldots, \alpha_n e_n)^t$ with $G = (g_{ij})_{d^*+1 \leq i \leq d; 1 \leq j \leq n} = \Psi_{21}\Psi_{11}^{-1}A^t_{(1)} - A^t_{(2)}$ that satisfy $GG^t = A^t_{(2)}(I - A_{(1)}\Psi_{11}^{-1}A^t_{(1)})A_{(2)}$. The matrix $B = I - A_{(1)}\Psi_{11}^{-1}A^t_{(1)}$ is symmetric and $B^2 = B$, then its eigenvalues are 0 and 1. Moreover, according to Lemma 2, on the event $\Omega_0$, the eigenvalues of $\Psi_{22} = A^t_{(2)}A_{(2)}$ are smaller than $6\mu_M$. Since, $GG^t = A^t_{(2)}BA_{(2)}$, the eigenvalues of $GG^t$ are smaller than $6\mu_M$. For $j \in \{d^*+1, \ldots, d\}$, this implies that $\sum_{k=1}^n g^2_{jk} = (GG^t)_{jj} \leq \sup_{u \in \mathbb{R}^{d-d^*}:\|u\|_2=1}\|GG^t u\|_2 \leq 6\mu_M$ and that $\zeta_j$ is a zero-mean Gaussian variable with variance satisfying

$$\mathbb{V}_{|\mathbb{X}}(\zeta_j) = \sigma^2 \sum_{k=1}^n g^2_{jk}\alpha^2_k \leq \frac{6\sigma^2 M_K \mu_M}{nh^d}, \tag{16}$$

where $\mathbb{V}_{|\mathbb{X}}$ stands for the variance symbol conditionally to $\mathbb{X} = (X_1, \ldots, X_n)$.

Moreover, for any $j \in \{d^* + 1, \ldots, d\}$, we have

$$|b_j| = \Big| \sum_{k=1}^{n} g_{jk} \Delta_k \Big| \leq \Big( \sum_{k=1}^{n} g_{jk}^2 \Big)^{1/2} \|\Delta\|_2 \leq \sqrt{6\mu_M} \|\Delta\|_2.$$

The $l_2$-norm of $\Delta$ can be upper bounded, with high probability, by using Bernstein's inequality. We have

$$\|\Delta\|_2^2 = \frac{1}{n} \sum_{i=1}^{n} V_i, \text{ where } V_i = \frac{1}{h^d} K\Big(\frac{X_i - x}{h}\Big) (f(X_i) - P_1 f(X_i, x))^2$$

and we use Assumption 5 to obtain $|V_i| \leq h^{2-d} L^2 M_K$ and $\mathbb{V}(V_i) \leq \mu_M h^{4-d} L^4 M_K$. Thus, there exists an event $\Omega_3$ of probability measure greater than $1 - \exp\big(-(3/8)\mu_M nh^d \mu_M^2 M_K^2\big)$, on which $\|\Delta\|_2^2 \leq \mathbb{E}\|\Delta\|_2^2 + \mu_M L^2 M_K h^2$. It is easy to see that $\mathbb{E}\|\Delta\|_2^2 \leq \mu_M L^2 M_K h^2$. Thus, on the event $\Omega_3$, we have $\|\Delta\|_2^2 \leq 2\mu_M L^2 M_K h^2$ and so $\max_{j=d^*+1,\ldots,d} |b_j| \leq 2\sqrt{3M_K}\mu_M Lh$.

We have $\lambda = 8\sqrt{3M_K}\mu_M Lh$ and, by using the classical upper bound on the tail of Gaussian random variables and (16), there exists a constant $c_0$ depending only on $\mu_M, M_K, L, \sigma$ such that

$$\mathbb{P}[\Omega_2^c | \mathbb{X} \in \Omega_0 \cap \Omega_3] \leq \sum_{j=d^*+1}^{d} \mathbb{P}[|\zeta_j| \geq \lambda/2 - |b_j| | \mathbb{X} \in \Omega_0 \cap \Omega_3]$$

$$\leq \sum_{j=d^*+1}^{d} \mathbb{P}[|\zeta_j| \geq \lambda/4 | \mathbb{X} \in \Omega_0 \cap \Omega_3]$$

$$\leq (d - d^*) \frac{c_0}{\sqrt{nh^{d+2}}} \exp\Big(-\frac{nh^{d+2}}{c_0^2}\Big),$$

where $\mathbb{P}[\cdot | \mathbb{X} \in \Omega_0 \cap \Omega_3]$ is the probability conditionally to $\mathbb{X}$ and to the event $\Omega_0 \cap \Omega_3$.

Here, we study the probability measure of $\Omega_1$. We have

$$\Omega_1^c \subset \cup_{j=0}^{d^*} \big\{ |\xi_j| \geq |\theta_j^*| - \lambda|a_j| - |b_j| \big\},$$

where $a = (a_0, \ldots, a_{d^*})^t = (\Psi_{11})^{-1} \overrightarrow{\text{sign}}(\theta_{(1)}^*), \xi = (\xi_0, \ldots, \xi_{d^*})^t = \Psi_{11}^{-1} A_{(1)}^t (\alpha_1 e_1, \ldots, \alpha_n e_n)^t$ and $b = (b_0, \ldots, b_{d^*})^t = \Psi_{11}^{-1} A_{(1)}^t (\Delta_1, \ldots, \Delta_n)^t$.

The random variable $\xi$ is of the form $H(\alpha_1 e_1, \ldots, \alpha_n e_n)^T$ with $H = (h_{ij})_{0 \leq i \leq d^*; 1 \leq j \leq n} = \Psi_{11}^{-1} A_{(1)}^t$ that satisfy $HH^t = \Psi_{11}^{-1}$. For $j = 0, \ldots, d^*$, the random variable $\xi_j$ is a zero-mean Gaussian variable with variance (conditionally to $X$) $\sigma_j^2 = \sigma^2 \sum_{k=1}^{n} h_{jk}^2 \alpha_k^2$ satisfying, on the event $\Omega_0$,

$$\sigma_j^2 \leq \frac{\sigma^2 M_K}{nh^d} (\Psi_{11}^{-1})_{jj} \leq \frac{8\sigma^2 M_K}{\mu_m nh^d}.$$

The last inequality holds, because, on the event $\Omega_0$, the maximum eigenvalue of $\Psi_{11}^{-1}$ is smaller than $8/\mu_m$.

Moreover, using Lemma 2, on the event $\Omega_0$, we have, for any $j \in \{0, \ldots, d^*\}$,

$$|a_j| \leq \|a\|_2 \leq \frac{8}{\mu_m} \|\overrightarrow{\mathrm{sign}}(\theta_{(1)}^*)\|_2 \leq \frac{8\sqrt{d^*+1}}{\mu_m} \leq \frac{8\sqrt{d_0+1}}{\mu_m}$$

and

$$|b_j| \leq \left( \sum_{k=1}^n h_{jk}^2 \right)^{1/2} \|\Delta\|_2 \leq (\Psi_{11}^{-1})_{jj} \|\Delta\|_2 \leq \frac{8}{\mu_m} \|\Delta\|_2.$$

We use the same argument as previously, to obtain that, on the event $\Omega_3$, the $l_2$-norm of $\Delta$ satisfies, $\|\Delta\|_2^2 \leq 2\mu_M L^2 M_K h^2$ and so $\max_{j=0,\ldots,d^*} |b_j| \leq (8/\mu_m)\sqrt{2\mu_M M_K} Lh$.

Since $\lambda = 8\sqrt{3M_K}\mu_M Lh$, we have $\lambda |a_j| + |b_j| \leq 36(\mu_M/\mu_m)LM_K\sqrt{d_0}h$. Thus, by using the classical upper bound on the tail of Gaussian random variables, the upper bound on the variance of the $\xi_j$'s and Assumption 6 with $C \geq 72(\mu_M/\mu_m)LM_K\sqrt{d_0}$, there exists a constant $c_1 > 0$ depending only on $\mu_m, \mu_M, M_K, C$ and $\sigma$ such that

$$\mathbb{P}[\Omega_1^c | \Omega_0 \cap \Omega_3] \leq \sum_{j=0}^{d^*} \mathbb{P}\left[ |\xi_j| \geq |\theta_j^*| - \lambda|a_j| - |b_j| \right] \leq \sum_{j=0}^{d^*} \mathbb{P}[|\xi_j| \geq Ch/2]$$

$$\leq 2(d^*+1)\frac{c_1}{\sqrt{nh^{d+2}}} \exp\left( -\frac{nh^{d+2}}{c_1^2} \right).$$

Thus, there exists a constant $c_2$ depending only on $\mu_m, \mu_M, M_K, L, C$ and $\sigma$ such that $\mathbb{P}[\Omega_1^c \cup \Omega_2^c | X \in \Omega_0 \cap \Omega_3] \leq d(c_2/\sqrt{nh^{d+2}}) \exp(-nh^{d+2}/c_2^2)$. Finally, using the results of Proposition 1 and Lemma 3, we obtain an upper bound on the probability measure of the event $\Omega_0$. Combining this upper bound and the result on the probability measure of $\Omega_3$, there exists $c_3$ and $c_4$ depending only on $L_\mu, \mu_m, \mu_M, M_K, L, C$ and $\sigma$ such that

$$\mathbb{P}[\Omega_0 \cap \Omega_1 \cap \Omega_2] \geq 1 - c_4 \exp(c_4 d - nh^{d+2}c_3).$$

$\square$

Theorem 1 follows by applying Proposition 2 and Proposition 3. For the study of $\tilde{\theta}$ when $f(x) = 0$, we just need to "move" the first coordinate of $\theta^*$ to the end. Namely, we can use the same arguments as previously for the following model

$$Z = \dot{A}\dot{\theta}^* + \varepsilon, \text{ where } \dot{\theta}^* = (\dot{\theta}_1^*, \ldots, \dot{\theta}_{d+1}^*)^t = (h\partial_1 f(x), \ldots, h\partial_d f(x), f(x))^t \quad (17)$$

and the lines of the design matrix $\dot{A} \in \mathcal{M}_{n,d+1}$ are $\dot{A}_i = \alpha_i \dot{U}\left(\frac{X_i - x}{h}\right), i = 1, \ldots, n$ (with $\dot{U}(v) = (v_1, \ldots, v_d, 1)$). In models (7) and (17), all the null coordinates of $\theta^*$ and $\dot{\theta}^*$ are at the end of the vector and all the non-zero coordinates are at the beginning of the vector. Then, we just have to consider $\dot{A} = (\dot{A}_{(1)}\dot{A}_{(2)})$, with $\dot{A}_{(1)} \in \mathcal{M}_{n,d^*}$ and $\dot{A}_{(2)} \in \mathcal{M}_{n,d-d^*+1}$, the matrices $\dot{\Psi} = \dot{A}^t\dot{A}$ and $\dot{\Psi}_{ij} = \dot{A}_{(i)}^t\dot{A}_{(j)}, i, j = 1, 2$. The proof, in this case, follows the line of the previous proof, with these notation and some minor changes in the indices.

### *4.2. Proof of Theorem 2*

Let $\delta > 0$. We have

$$\mathbb{P}[|\widehat{f}(x) - f(x)| \geq \delta]$$

$$= \mathbb{P}\Big[|\widehat{f}(x) - f(x)| \geq \delta | \widehat{J}_2 = J\Big] \mathbb{P}[\widehat{J}_2 = J] + \mathbb{P}\Big[|\widehat{f}(x) - f(x)| \geq \delta | \widehat{J}_2 \neq J\Big] \mathbb{P}[\widehat{J}_2 \neq J]$$

$$\leq \mathbb{P}\Big[|\widehat{f}(x) - f(x)| \geq \delta | \widehat{J}_2 = J\Big] + \mathbb{P}[\widehat{J}_2 \neq J] \mathbb{1}_{\delta \leq (2f_{max})^2}$$

$$\leq c_1 \exp(-c_2 n^{\frac{2\beta}{2\beta + d^*}} \delta^2) + c_1 \exp(c_1 d - nh^{d+2} c_0) \mathbb{1}_{\delta \leq (2f_{max})^2},$$

where, on the event $\{\widehat{J}_2 = J\}$, we used the classical result on LPE (cf. Audibert and Tsybakov (2007)) and, for the event $\{\widehat{J}_2 \neq J\}$, we upper bounded its probability measure by using Theorem 1. The assumption on $d$ completes the proof. $\qquad\square$

### Acknowledgment

Karine Bertin is supported by Project FONDECYT 1061184, project CNRS CONICYT CNRS0705 and project PBCT ACT 13 laboratory ANESTOC. The authors would like to thank Alexander Tsybakov for his contribution to this work.